\documentclass[conference,a4paper]{IEEEtran}

\pdfoutput=1

\usepackage{algorithm}
\usepackage{algorithmic}
\usepackage{amsmath,amssymb,amsfonts}
\usepackage{graphicx}
\usepackage[usenames]{color}

\newtheorem{proposition}{Proposition}
\newtheorem{theorem}{Theorem}

\newtheorem{remark}{Remark}

\def\EE{{\mathbb E}}

\def\Ex#1{{\EE\left\{#1\right\}}}

\def\CC{\mathbb C}
\def\ZZ{\mathbb Z}

\def\be{\begin{equation}}
\def\beq#1{\begin{equation}\label{#1}}
\def\ee{\end{equation}}
\def\bea{\begin{eqnarray}}
\def\beqa#1{\begin{eqnarray}\label{#1}}
\def\eea{\end{eqnarray}}
\def\ba{\begin{array}}
\def\ea{\end{array}}

\DeclareMathAlphabet{\mathpzc}{OT1}{pzc}{m}{it}

\def\cC{{\mathcal C}}
\def\cCN{{\mathcal {C\!N}}}

\def\cG{{\mathcal G}}

\def\cL{{\mathcal L}}

\def\cN{{\mathcal N}}

\def\cS{{\mathcal S}}

\def\cV{{\mathcal V}}

\def\ba{{\mathbf a}}

\def\bg{{\mathbf g}}
\def\bG{{\mathbf G}}

\def\bN{{\mathbf N}}

\def\bU{{\mathbf U}}

\def\bx{{\mathbf x}}
\def\bX{{\mathbf X}}

\def\bY{{\mathbf Y}}
\def\bZ{{\mathbf Z}}

\renewcommand\bg{{\boldsymbol g}}
\renewcommand\bx{{\boldsymbol x}}

\renewcommand\bN{{\boldsymbol N}}
\renewcommand\bU{{\boldsymbol U}}
\renewcommand\bX{{\boldsymbol X}}
\renewcommand\bY{{\boldsymbol Y}}
\renewcommand\bZ{{\boldsymbol Z}}
\renewcommand\Re{\mathsf{Re}}

\def\cbG{\boldsymbol\cG}

\DeclareMathOperator*{\argmin}{arg\,min}

\DeclareGraphicsExtensions{.pdf,.jpeg,.png}

\begin{document}
%
\title{Estimation of frequency modulations on wideband signals; applications to
  audio signal analysis}

\author{\IEEEauthorblockN{Harold Omer and Bruno Torr\'esani}
\IEEEauthorblockA{Aix-Marseille Universit\'e, CNRS, Centrale Marseille, LATP,
  UMR 7353, 13453 Marseille, France\\
Email: harold.omer@latp.univ-mrs.fr\quad
bruno.torresani@latp.univ-mrs.fr}
}

\maketitle

\begin{abstract}
The problem of joint estimation of power spectrum and modulation from
realizations of frequency modulated stationary wideband signals is considered.
The study is motivated by some specific signal classes from which departures to
stationarity can carry relevant information and has to be estimated.

The estimation procedure is based upon explicit modeling of the signal as a
wideband stationary Gaussian signal, transformed by time-dependent, smooth
frequency modulation. Under such assumptions, an approximate expression for the
second order statistics of the transformed signal's Gabor transform is obtained,
which leads to an approximate maximum likelihood estimation procedure. 

The proposed approach is validated on numerical simulations.
\end{abstract}


%
\IEEEpeerreviewmaketitle

\section{Introduction}
Usual time-frequency models for audio signals often rest upon expansions with
respect to dictionaries of time-frequency waveforms, such as Gabor frames,
wavelet frames, or more general families. Such descriptions are generally
adequate for signal classes such as (voiced) speech, music,... where
specific time-frequency localisation properties can be exploited. They are less
effective for less structured signals, such as wideband sound signals.

We are concerned here with an alternative description of audio signals, aiming
at describing different sound classes such as environmental noise, engine
sound,... which are in addition non-stationary, in the sense that they carry
information related to dynamics. As an example, think of an accelerating engine
sound, where the acceleration can generally be perceived. This example suggests
to study sound models, which we will term {\em timbre$\times$dynamics}, in
which a reference (stochastic) stationary signal, characterized by its timbre,
is modulated by some dynamic deformation. Given such signals, a problem is to
estimate the modulation (and possibly the underlying power spectrum). While many
techniques have been developed for frequency modulation estimation for narrow
band signals (see e.g.~\cite{VanTrees03detection}, the wideband case is more
complex and has apparently received less attention.

A class of models based upon deformations of stationary processes has been
proposed and studied in~\cite{Clerc03estimating}, motivated by the famous {\em
  shape from texture} image processing problem. A main aspect of the approach is
based on the remark that a generic class of transformations can be represented
by transport equations in a suitable representation space.

We adopt here a more explicit point of view, and limit to stationary Gaussian
processes, transformed by a time-dependent modulation. We characterize the
distribution of fixed time slices of a Gabor transform of such signals, and
formulate the corresponding maximum likelihood estimation problem. As a result,
we provide an estimation algorithm which is demonstrated on a small number of
numerical examples.

\section{Frequency modulation of stationary random signals}

\subsection{Notations and background}
\subsubsection{Random signals}
We shall be concerned with complex Gaussian random signal models $\bX$ of finite
length $L$, which we shall assume zero-mean for the sake of simplicity. As is
customary in finite-dimensional Gabor analysis, we shall
also assume periodic boundary conditions, i.e. $\bX_{t+L}=\bX_t$.
Given such a signal $\bX$, we shall denote by $C_{\bX}$ its covariance matrix, and by
$R_\bX$ its relation matrix (see~\cite{Picinbono96second} for details), defined as
\beq{fo:cov.rel.mat}
C_\bX(t,s) = \Ex{X_t\overline{X}_s}\ ,\quad
R_\bX(t,s) = \Ex{X_tX_s}\ ,
\ee
and we will write $\bX\sim\cCN(0,C_\bX,R_\bX)$. $\bX$ is said to be {\em circular}
if $R_\bX=0$.
\subsubsection{Time-frequency representation}
We shall use the following notations. Given a window function $\bg$, the
corresponding short time Fourier transform of a signal (STFT) $\bx\in\CC^L$ is defined by
\beq{fo:STFT}
\cV_\bg\bx(m,n) = \sum_{t=0}^{L-1} x[t] \overline{g}[t-n]e^{-2i\pi m(t-n)/L}\ .
\ee
Given lattice constants $a$ and $b$ (divisors of the signal length $L$)
the corresponding Gabor transform reads
\beq{fo:DGT}
\cG_\bx[m,n] = \cV_\bg \bx(mb,na)\ m=0,\dots M-1,\ n=0\dots N-1\ ,
\ee
with $M=L/b$ and $N=L/a$. $\cG_\bX$ is an $M\times N$ array.
For suitably chosen $\bg$, and $a$ and $b$ small
enough, the Gabor transform is invertible
(see~\cite{Carmona98practical,Grochenig01foundations});
in finite dimensional situations, efficient
algorithms have been developed and implemented
(see~\cite{ltfatnote011}).

\begin{remark}[Notations]
As usual, summation bounds in the frequency domain depend of the parity of the
signal length $L$. For the sake of simplicity, we introduce some notations and
denote by $I_L$ the integer interval $I_L=[1-L/2,L/2]$ if $L$ is even, and the
integer interval $I_L=[-(L-1)/2, (L-1)/2]$ if $L$ is odd.
The corresponding positive frequencies interval will be denoted by
$I_L^+=[0,L/2]$ if $L$ is even, and $I_L^+=[0,(L-1)/2]$ if $L$ is odd.
\end{remark}

\subsection{The model: definition and main estimates}
We are concerned here in a simple model of signal transformation, which may be
written as follows. We denote by $\bX$ a zero-mean, wide sense stationary
Gaussian random process, with covariance matrix $C_\bX$, and by
$\bZ$ the associated analytic signal.
We denote by $\cS_\bX$ the power spectrum of
$\bX$, and assume that $\cS_\bX(0)=0$, and if $L$ is even, that $\cS_\bX(L/2)=0$.
Under such an assumption, it is easy to show that $\bZ$ is a circular complex
Gaussian random vector (by a finite dimensional version of a standard argument,
see e.g.~\cite{Picinbono94circularity}).

The observation is assumed to be the real part $\bY_r = \Re(\bY)$ of a complex
valued signal $\bY$; for the sake of simplicity we shall only work with the
latter, assumed to be an USB (upper sideband) modulated version $\bY$ of a
reference stationary signal $\bX$, of the form
\begin{equation}
\label{fo:signal.model}
Y_t = Z_t e^{2i\pi\gamma(t)/L} + N_t\ ,
\end{equation}
where $\gamma\in C^2$ is an unknown smooth, slowly varying {\em modulation
  function}, and $\bN=\{N_t,\,t=0,\dots L-1\}$ is a real
Gaussian white noise, with variance $\sigma_0^2$.
Obviously, when $\gamma$ is not a constant function, $\bY$ is not a wide sense
stationary signal any more. The problem at hand is to estimate the unknown
modulation $\gamma$ and the original power spectrum $\cS_\bX$ from a single
realization of $\bY$.

\medskip
Clearly, $\bZ\sim\cC\cN(0,C_\bZ,0)$ is a circular complex Gaussian
random signal, with covariance matrix
\begin{equation}
C_\bZ(t,s) = \sum_{\nu \in I_L^+} \cS_\bX(\nu) e^{2i\pi\nu(t-s)/L}\ ,
\end{equation}
and is therefore wide-sense stationary.

In the proposed approach, we will base the estimation on a Gabor representation
of the observed signal, and deliberately disregard correlations across time of
the Gabor transform (hence focusing on time slices of the Gabor transform of
the observation).
The distribution of time slices of the analytic signal $\bZ$ of the
original signal is characterized in the following two results, which result from
direct calculations.
\begin{proposition}
For fixed $n$, the Gabor transform $\cG_{\bN}[.,n]$ of the gaussian white noise
is a stationary Gaussian random vector, with circular covariance matrix
\begin{equation}
C_{\cG_\bN}[m,m'] = \sigma_0^2\sum_{k=0}^{L-1} \overline{\hat g}[k]\hat g[k-(m'-m)b]
\end{equation}
\end{proposition}
\begin{proposition}
For fixed time index $n$, the Gabor transform $\cG_{\bZ}[.,n]$ of the analytic
signal is a circular complex Gaussian random vector, with covariance matrix
\begin{equation}
C_{\cG_\bZ}[m,m'] = \sum_{k \in I_L^+} \cS_\bX[k] \overline{\hat g}[k-mb]\hat g[k-m'b]
\end{equation}
\end{proposition}
The estimation of the modulation will be based upon an approximation of the
covariance matrix of the observed signal. In a few words, the Gabor transform of
the frequency modulated signal can be approximated by a deformed version of the
Gabor transform of the original signal. The deformation takes the form of a
time-varying frequency shift. A more precise argument, based upon first
order approximation of the modulation function $\gamma$, leads to the
following result.
\begin{theorem}
\begin{enumerate}
\item
For fixed time, the Gabor transform $\cG_\bY$ may be approximated as
\begin{equation}
\cG_\bY[m,n]= \bG^{(n;\gamma'(na)/b)}[m] + R[m]\ \label{eq:approximation1},
\end{equation}
where $\bG^{(n;\delta)}$ is a frequency-shifted Gabor transform
\begin{eqnarray}
\bG^{(n;\delta)}[m] = \sum_{t=0}^{L-1} Z_t
\overline{g}[t-na] e^{-2i\pi[m-\delta][t-an]/M } \nonumber \\
 + \ \cG_\bN[m,n] \ , \label{eq:approximation2}
\end{eqnarray}
and the remainder is bounded as follows: for all $m,m'$,
\begin{equation}
\left|\Ex{R[m]\overline{R}[m']}\right| \le
\sigma^2_Z \left( \frac{\pi e}{L} \| \gamma''\|_{\infty} \mu_2 +
2\mu_1 \right)^2\!\!,
\end{equation}
where  $\sigma_Z^2$ is the variance of $Z$ and with 
\begin{equation}
\mu_1 = \sum_{t\in I_T^c}|g(t)| \ , \; \mu_2 = \sum_{t \in I_T} t^2 |g(t)|\ ,\;
T = \sqrt{\frac{L}{\pi {\| \gamma'' \|}_{\infty}}}
\end{equation}
where $I_T=[-T,T]$ and $I_T^c=I_L\backslash I_T$

\item
Given $\delta$, and for fixed $n$, $\bG^{(n;\delta)}$ is distributed
following a circular multivariate complex Gaussian law, with covariance matrix
\begin{equation}
C_{\bG^{(n;\delta)}}[m,m']=C_{\cG_\bZ}[m-\delta,m'-\delta]
+ \ C_{\cG_\bN}[m,m']\ .
\label{eq:cov}
\end{equation}
\end{enumerate}
\end{theorem}
The estimation procedure described below is a maximum likelihood
approach, which requires inverting the covariance matrix of
vectors $\bG^{(n;\delta)}$.
The latter is positive semi-definite by construction, but not necessarily
definite. The result below provides a
sufficient condition on $\bg$ and the noise for invertibility.
\begin{proposition}
Assume that the window $\bg$ is such that
\begin{equation}
K_\bg:=\min_{t=0\dots L-1}\left(\sum_{k=0}^{b-1}|g[t+kM]|^2\right)>0\ .
\end{equation}
Then for all $\bx\in\CC^M$,
\begin{equation}
\bx^* C_{\bG}\bx \ge \sigma_0^2 K_\bg\ ,
\end{equation}
and the covariance matrix is therefore boundedly invertible.
\end{proposition}
\begin{remark}
\label{rem:condition}
The condition may seem at first sight unnatural to Gabor
frame experts. However, it simply expresses that the number $M$ of
frequency bins shouldn't be too large
if one wants the covariance matrix to be invertible. However, reducing $M$ also
reduces the precision of the estimate, and a trade-off has to be found, as
discussed in the next section.
\end{remark}

\subsection{Improving the frequency resolution}
\label{subsec:improving}
We propose here a method to improve the frequency resolution of our estimations.
We have already seen that the invertibility of the covariance matrix 
requires that the number of frequency bins of the Gabor transform shouldn't be
too large. As a result
however, it may be convenient, as we shall see later, to have access to 
the information contained in all the frequency frames of the short time 
Fourier transform defined in equation (\ref{fo:STFT}). For this purpose, 
we also consider alternative versions of the Gabor transform, associated with
frequency-shifted sampling lattices:
\beq{fo:DGT2}
\cG^c_\bx[m,n] = \cV_\bg \bx(mb+c,na)\ ,\quad
m\in\ZZ_M,\ n\in\ZZ_N\ ,
\ee
where $c \in [0 , b-1]$ . We now have at our disposal a collection of $b$ Gabor
transforms, which are all different subsampled versions of the STFT.
The previous results and proofs remain valid with this new definition of the
Gabor transform. Equations~\eqref{eq:approximation1}
and~\eqref{eq:approximation2} now become
\begin{equation}
\cG^c_\bY[m,n]= \bG^{(n;\gamma'(na)/b+c/b)}[m] + R\ ,
\end{equation}
where
\begin{eqnarray}
\bG^{(n;\delta^c)}[m] &=& \sum_{t=0}^{L-1} Z_t
\overline{g}[t-na] e^{-2i\pi[m-\delta^c][t-an]/M } \nonumber \\
&&\hphantom{aaaaa} +\ \cG^c_\bN[m,n] \ , 
\end{eqnarray}
and the associated Equation~\eqref{eq:cov} now reads:
\be
C_{\bG^{(n;\delta^c)}}[m,m']=C_{\cG_\bZ}[m-\delta^c,m'-\delta^c]
+ \ C_{\cG_\bN}[m,m']\, .
\ee

The rationale will be that a frequency shift $\delta^c$ can be estimated from
each one of these thansforms, and the optimal one will be retained.

\section{Estimation procedure}
We now describe in some details the estimation procedure corresponding
to our problem. The estimation problem is the following: from
a single realization of the signal model~\eqref{fo:signal.model}, estimate
the modulation function $\gamma$ and the original power spectrum
$\cS_\bX$. We first notice the indeterminacy in the problem, namely the
fact that adding an affine function to $\gamma$ is equivalent to
shifting $\cS_\bX$. This has to be fixed by adding an extra constraint
in the estimation procedure.

\subsection{Maximum likelihood modulation estimation}
We now turn to the estimation procedure, that exploits the above results.
With the same notations as before, we fix a value of the time index $n$, and
denote for simplicity by $\cbG=\cbG^{(n)}$ the corresponding fixed time slice of
$\cG^c_\bZ$.
Due to the multivariate complex Gaussian distribution of the signal and the
fixed time Gabor transform slices, the log-likelihood of a slice takes the form
\begin{equation}
\cL_\delta(\cbG) = \cbG^*\left(C_{\bG^{(n;\delta^c)}}\right)^{-1} \cbG +
\ln\left(\pi^M \det(C_{\bG^{(n;\delta^c)}})\right)\ .
\end{equation}
Therefore, the maximum likelihood estimate for the frequency shift assumes the
form
\begin{equation}
\hat\delta^c\! =\!
\argmin_{\delta^c}\left[\cbG^*\!\left(C_{\bG^{(n;\delta^c)}}\right)^{-1}\! \cbG + 
\ln\left(\!\pi^M\!\! \det(\!C_{\bG^{(n;\delta^c)}}\!)\!\right)\!\right]\, .
\end{equation}
However, we notice that $\det(C_{\bG^{(n;\delta^c)}})$ actually does not depend on
the modulation parameter $\delta^c$. Therefore the maximum likelihood estimate
reduces to
\begin{equation}
\hat\delta^c = \argmin_{\delta^c} \left[\cbG^*\left(C_{\bG^{(n;\delta^c)}}\right)^{-1}
  \cbG\right]\ , \label{eq:ml1}
\end{equation}
a problem to be solved numerically. Notice that this requires the knowledge of
the covariance matrix $C_{\bG^{(n;0)}}$ corresponding to the Gabor transform of
the noisy stationary signal. The latter is generally not available, and has to
be estimated as well.

As $\delta^c(n)\approx \left( \gamma'(an) - c \right)/b$, the estimates of
$\delta$ for each $n$ lead to an estimate of $\gamma'$.
Since we solve the minimisation problem by an exhaustive search on the $\delta^c$, 
the estimate of $\gamma'$ is coarsely quantized (see Remark~\ref{rem:condition}),
as $b$ is large and $\hat\gamma'(an) \in [c, b+c,2b+c,..,(M-1)b+c]$. 
This problem is solved by using the family of frequency-shifted versions of
Gabor transform described in subsection~\ref{subsec:improving} and making a new
exhaustive search on the $\hat\delta^c$
\begin{equation}
\hat\delta = \argmin_{c} \left[\cbG^*\left(C_{\bG^{(n;\hat\delta^c)}}\right)^{-1}
  \cbG\right]\ . \label{eq:ml2}
\end{equation}
The quantization effect on the final estimation of the modulation function
is therefore attenuated, i.e. $\hat\gamma'(an) \in [0, L-1]$.
Obtaining from this estimation a smoother estimate for the 
modulation function $\gamma$ requires extra interpolation techniques.
\begin{remark}
As an alternative, one may also avoid exhaustive searches and seek minimizers
in~\eqref{eq:ml1} using more elaborate numerical techniques, that would avoid
quantization effects. This question is currently under study.
\end{remark}

\subsection{Estimation of the underlying covariance matrix}
We now describe a method for estimating the covariance matrix $C_{\bG^{(n;0)}}$. 
Suppose that an estimate $\hat\gamma$ of the modulation function $\gamma$
is available. Then the signal $\bY$ can be demodulated by setting
\begin{equation}
\bU = \bY e^{-2i\pi\hat\gamma/L}\ , \label{eq:demodulation}
\end{equation}
Clearly, $\bU$ is an estimator of $\bZ  + \bN e^{-2i\pi\gamma/L}$, the noisy
stationary signal.
We can now compute the covariance matrix $C_{\cG_\bU}$ of the Gabor transform of
$\bU$, which is an estimator of $C_{\cG_\bZ} + C_{\cG_\bN}$. 
Comparing with equation~(\ref{eq:approximation2}) we finally obtain an 
estimator for the covariance matrix
\begin{equation}
C_{\cG_\bU} \approx C_{\bG^{(n;0)}}
\end{equation}

\begin{remark}
The power spectrum $\cS_\bX$ of the stationary signal can be estimated from
$\bU$ using a standard Welch periodogram estimator, or by marginalizing the
square modulus of the Gabor transform of the demodulated signal, as
described in~\cite{Carmona98practical}.
\end{remark}


\subsection{Summary of the estimation procedure}
We now summarize an iterative algorithm to jointly estimate the covariance
matrix $C_{\bG^{(n;0)}}$ and the modulation function $\gamma$, that exploits
alternatively the two procedures described above.
The procedure is as follows, given a first estimation of the modulation
function, we can perform a first estimation of the covariance matrix, which in
turn allows us obtain a new estimation of the modulation function.
The operation is repeated until the stopping criterion is satisfied.

For the initialization, we need a first modulation frequency estimate, for
which we use the center of mass of the modulated signal Gabor transform
\beq{fo:initialization}
\hat{\delta}^{(0)}(n)= \frac{\sum_{m=0}^{M-1}
  m|\bG^{(n;\delta)}|^2[m]}{\sum_{m=0}^{M-1} |\bG^{(n;\delta)}|^2[m]} \ .
\ee
The stopping criterion is based upon the evolution of the frequency modulation
along the iterations. More precisely, we use the empirical criterion
\begin{equation}
\frac{||\hat\delta^{(k)} -\hat\delta^{(k+1)}||_2}{|| \hat\delta^{(k+1)}||_2} <
\epsilon \label{eq:criterion}
\end{equation}
The pseudo-code of the algorithm can be found below
\begin{algorithm}
\caption{Joint covariance and modulation estimation}
\label{alg:estimation}
   
\begin{algorithmic}

\STATE {\bf Initialize} as in~\eqref{fo:initialization}
\WHILE{ criterion (\ref{eq:criterion}) is false}
\STATE $\bullet$ Compute $\hat\gamma^{(k)}$ by interpolation from $\delta^{(k)}$.

\STATE $\bullet$ Demodulate $\bY$ using $\hat\gamma^{(k)}$
following~\eqref{eq:demodulation}

\STATE $\bullet$ Compute the Gabor transform of the demodulated signal
${\hat\bG}^{(k;n)}[m] = \cG_{{\bU}^{(k)}}[m ,n ]$

\STATE $\bullet$ Estimate $\hat{\delta}^{(k+1)}$ using the covariance matrix of ${\hat\bG}^{(k;n)}$ from (\ref{eq:ml1}) and (\ref{eq:ml2})

\STATE $\bullet$ $k:= k+1$

\ENDWHILE

\end{algorithmic}
\end{algorithm}

\section{Numerical results}
The  proposed estimation procedure has been implemented using
{\sc Matlab/Octave}, and relies on the {\sc Ltfat}
toolbox~\cite{Soendergaard12linear} for the time-frequency transforms.

\begin{center}
\begin{figure}
\includegraphics[width=8.5cm]{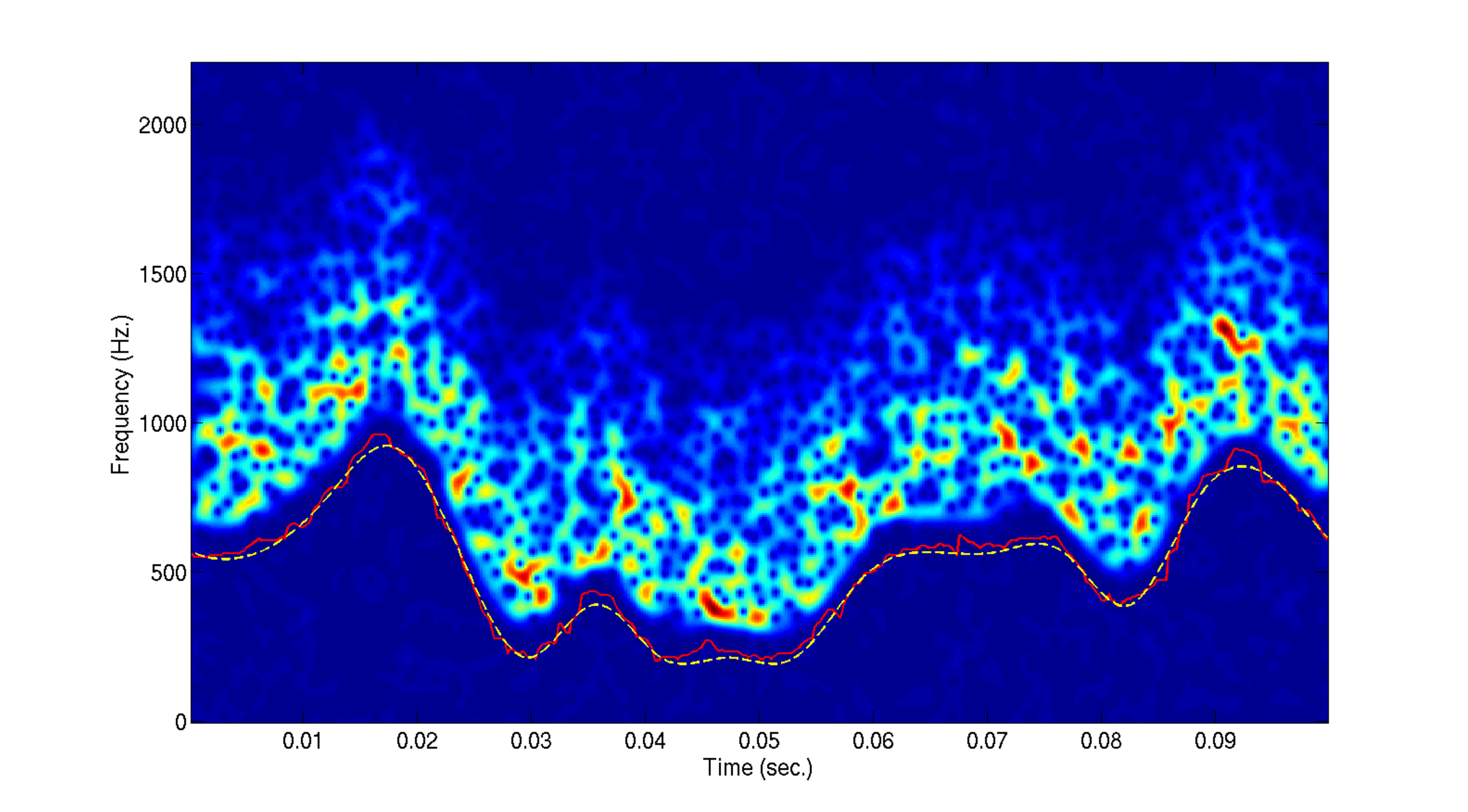}

\caption{Gabor transform of a frequency modulated synthetic stationary random
  signal, superimposed with the frequency modulation:  estimate (red) and
  original (yellow). \label{fig:estimation}}
\end{figure}
\end{center}

We display in Fig.~\ref{fig:estimation} an example of estimation result. The
original signal was generated as pseudo-random stationary Gaussian signal with a
smooth, wideband power spectrum, that was further modulated by a smooth
frequency modulation function.  Fig.~\ref{fig:estimation} displays the Gabor
transform of the modulated signal (positive frequencies only), together with the
original and the estimate for the frequency modulation. For the sake of clarity,
the frequency estimate has been displayed below the relevant part of the Gabor
transform (remember that it is defined up to an additive constant).
As can be seen, the result is fairly satisfactory, the estimated modulation
follows closely the ground truth.

\begin{center}
\begin{figure}
\centerline{\includegraphics[width=6cm]{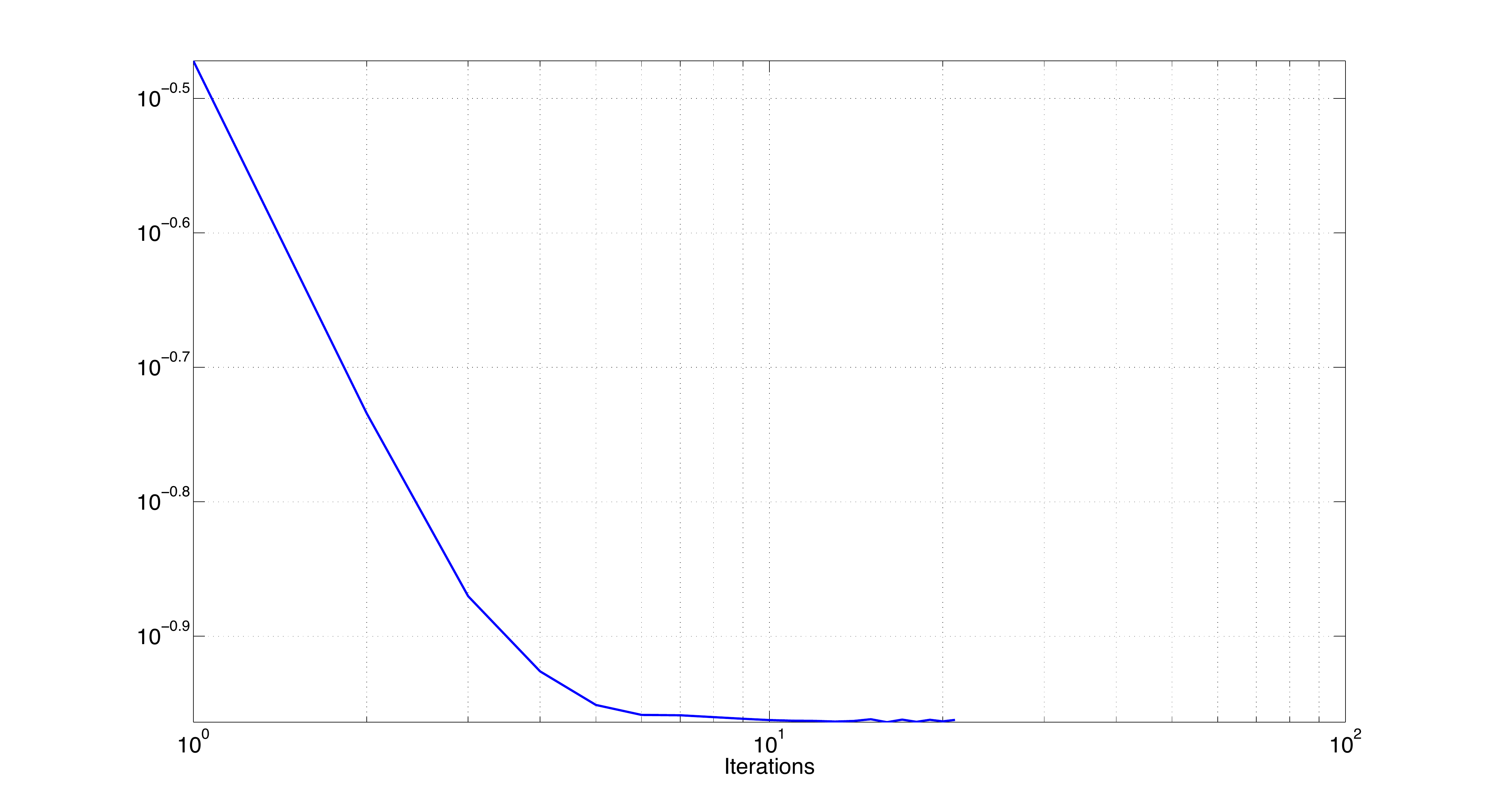}}
\caption{Log-log plot of the evolution of the criterion proposed in (\ref{eq:criterion}) according on the number of iterations. \label{fig:convergence} }
\end{figure}
\end{center}

To asses the convergence properties of the proposed approach, the same
experiment was run several times with the same modulation law and different
seeds for the underlying stationary noise. We display in
Fig.~\ref{fig:convergence} the evolution of the criterion as a function of the
iteration index, averaged over 20 realizations. Convergence appears to be fast,
with power-law like decay speed.

\section{Conclusion}
We have presented in this paper a new approach for modulation frequency and
power spectrum estimation from wideband signals, based upon explicit modeling.
A main point that is exploited in our approach is the fact that modulations can
be locally approximated by frequency shifts in the Gabor domain.
The algorithm has been validated using numerical simulations, that show that
when signals are generated according to the model of interest, very accurate
results can be obtained.

Further developments include numerical tests on real signals, such as natural
sounds generated by rolling bodies with variable speed,...
We shall also consider extending this approach to other transformation models,
such as time warping or more general transformations.


\section*{Acknowledgment}

This work was supported by the ANR project Metason ANR-10-CORD-010.



\bibliographystyle{IEEEtran}
\bibliography{IEEEabrv,OT}
%

%

\end{document}